\documentclass{article}
\usepackage[affil-it]{authblk}
\usepackage[english]{babel}
\usepackage[T2A]{fontenc}
\usepackage[cp1251]{inputenc}
\usepackage[tbtags]{amsmath}
\usepackage{amsfonts,amssymb,mathrsfs,amscd, graphicx}
\def\R{{\mathbb R}}
\newcommand{\SO}{\operatorname{SO}}
\newcommand{\SE}{\operatorname{SE}}
\newcommand{\spann}{\operatorname{span}}

\newtheorem{remark}{Remark}
\newtheorem{theorem}{Theorem}

\begin{document}

\title{Liouville integrability in a four-dimensional model of the visual cortex 
%\footnote{ The reported study of A.P. Mashtakov (chapters 1, 2, 4.2, 5) was funded by RFBR, project number 19-31-51023 and was carried out at the Sirius University of Science and Technology (Sochi, Russia) and Ailamazyan Program Systems Institute of RAS (Pereslavl-Zalessky, Russia)}
}
\author{I.A. Galyaev$^{1,2,*}$,
A.P. Mashtakov$^{1,3,**}$}
\affil{$^1$ Sirius University of Science and Technology\\
 Sochi, Adler district, Olympic Avenue, 1,\\
$^2$~V. A. Trapeznikov Institute of Control Sciences of RAS\\
Moscow, st. Profsoyuznaya, 65\\
$^3$~Ailamazyan Program Systems Institute of RAS,\\
 Yaroslavl region, Pereslavsky district, Veskovo, st.~Petera I, 4a\\
%Научно--технологический университет ``Сириус'',\\
%г. Сочи, Адлерский район, Олимпийский проспект, д. 1,\\
$^*$E-mail: ivan.galyaev@yandex.ru\\
$^{**}$E-mail: alexey.mashtakov@gmail.com}
\maketitle
\begin{abstract}
We consider a natural extension of the Petitot-Citti-Sarti model of the primary visual cortex. In the extended model, curvature of contours is taking into account. The occluded contours are completed via sub-Riemannian geodesics in the four-dimensional space M of positions, orientations, and curvatures. Here, $M = \R^2 \times \SO(2) \times \R$ models the configuration space of neurons of the visual cortex. We study the problem of sub-Riemannian geodesics on $M$ via methods of geometric control theory. We prove complete controllability of the system and existence of optimal controls. By application of Pontryagin maximum principle, we derive a Hamiltonian system that describes the geodesics. We obtain explicit parametrization of abnormal extremals. In the normal case, we provide three functionally independent first integrals. Numerical simulations indicate existence of one more first integral that results in Liouville integrability of the system. 
\end{abstract}

\section{Introduction}
\begin{figure}[ht]
\centering
\begin{minipage}{0.25\linewidth}
\includegraphics[width=\linewidth]{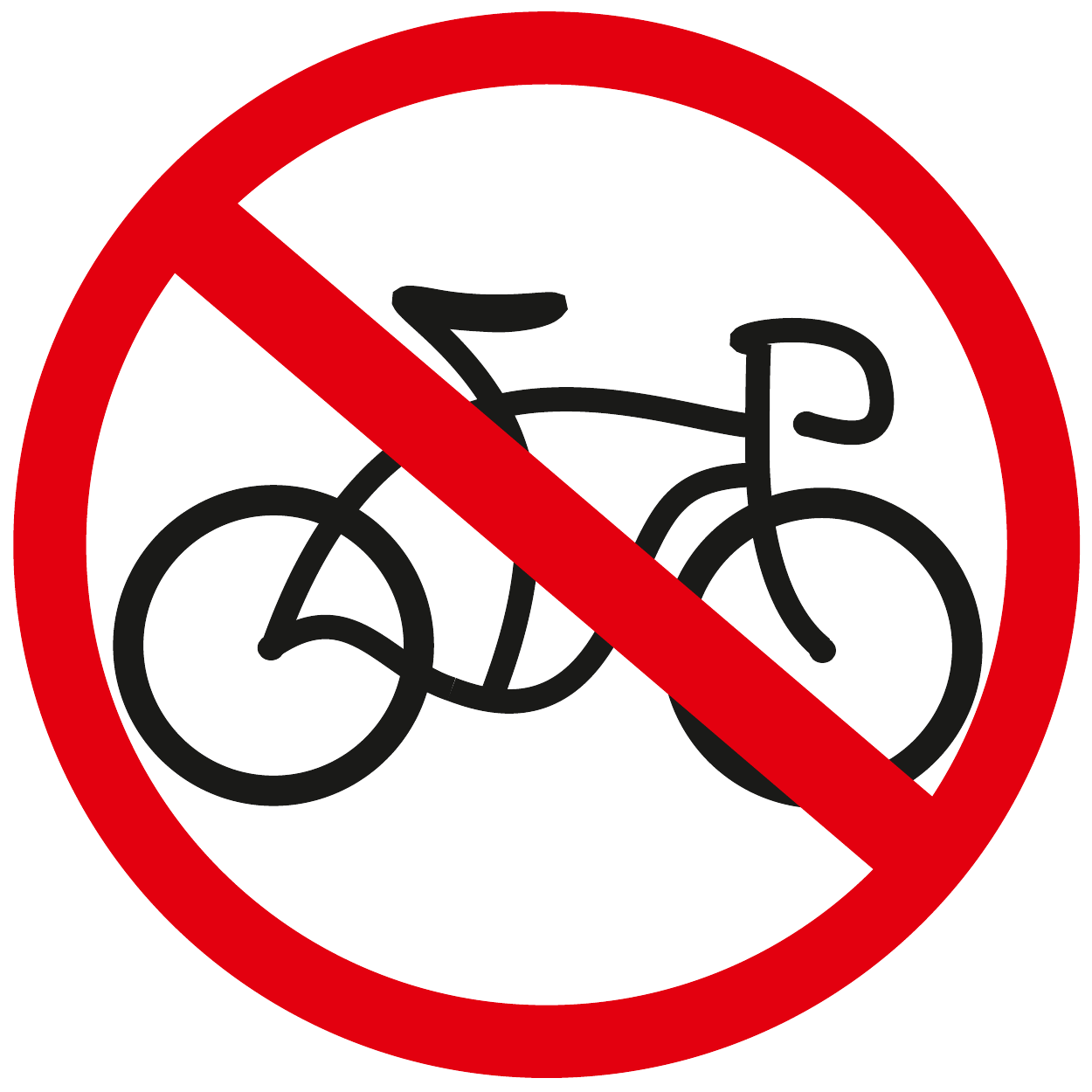}
\includegraphics[width=\linewidth]{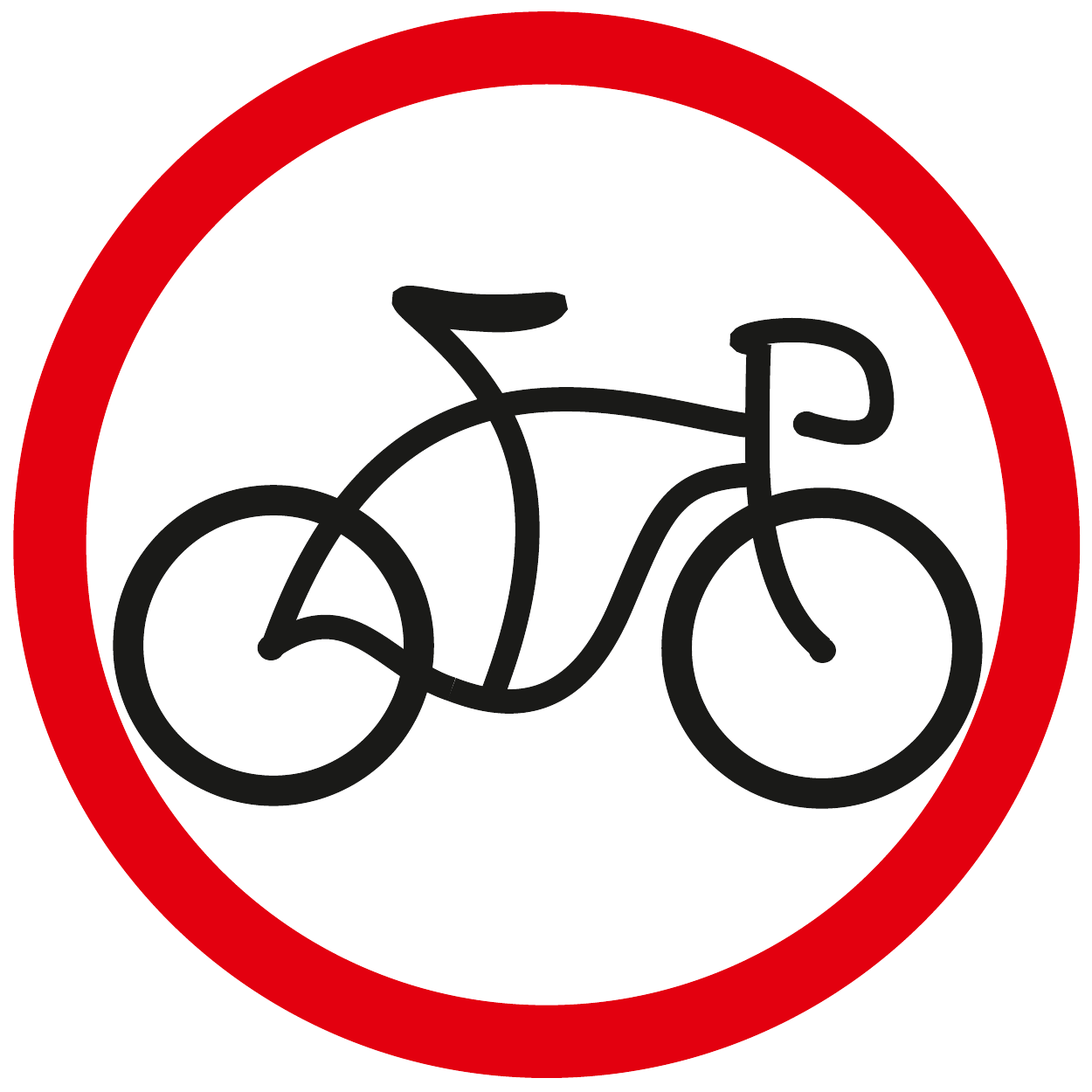}
\end{minipage}~
\begin{minipage}{0.74\linewidth}
\includegraphics[width=\linewidth]{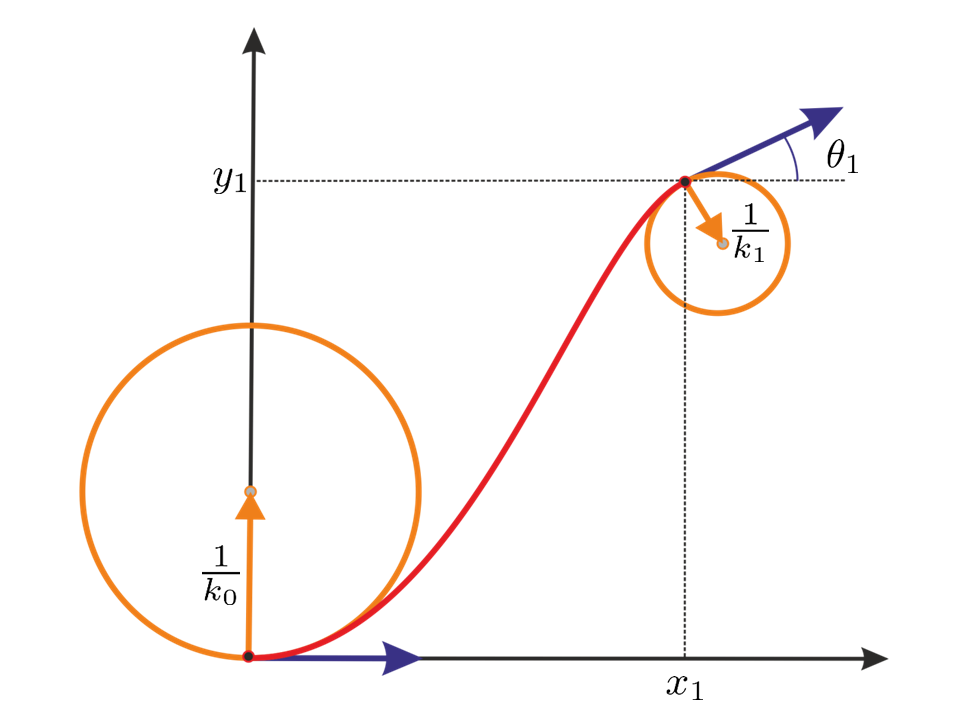}
\end{minipage}
\caption{In the four-dimensional model of the visual cortex an occluded contour is completed via the planar projection of a sub-Riemannian length-minimizer in the space $\R^2\times SO(2) \times \R = M\ni(x,y,\theta,k)$ of positions, orientations and curvatures. In the left column we show an example of the image with partially occluded contours and the complete image. In the right column we show a trajectory that satisfies given boundary conditions. The curvature is visualized as its reciprocal --- the radius of the osculating circle.}
\label{fig:Stat}
\end{figure}

The principles of biological visual systems are of great interest among researchers in many fields of science. One of the important directions is development and study of realistic mathematical models describing a certain stage of visual signal processing. In this paper, we consider the mathematical model of the visual cortex that was introduced and developed in~\cite{petit, petitot, ccc}. Such a model is obtained by extension of the classical Petitot-Citti-Sarti model, which represents the primary visual cortex V1 as a sub-Riemannian structure on the Lie group $\SE(2)=\R^2\times\SO(2)$. Here $ \SE(2) $ models the configuration space of neurons V1, which can be understood as the space of positions $ \R^2 $ and orientations $ \SO(2) $. The extension is performed by taking into account the curvature of contours of observed image. This leads to a sub-Riemannian structure in the 4-manifold  $M=\R^2 \times \SO(2) \times \R$, where the fourth component means the curvature. In the four-dimensional model of the visual cortex a partially occluded contour is completed via the planar projection of a sub-Riemannian length-minimizer in $M$ that satisfies the boundary conditions (position, orientation and curvature) obtained from the boundary of occluded region. See Fig.~\ref{fig:Stat}.

Interest in such models also comes from applications. The principles of biological visual systems are actively used in computer vision. Based on these principles, effective methods of image processing are created: enhancement, segmentation, inpainting, feature detection. For example, in~\cite{Duits2007}, the aurhors describe an approach that is based on lifting an image into the expanded space of positions and orientations. After such lifting, the sub-Riemannian length minimizers are used to detect the salient lines~\cite{SIAM, SO3JMIV, SE3JDCS}. This approach is actively used in medical image processing.

The mathematical model of V1 as a sub-Riemannian structure in the space of positions and orientations was proposed by J.~Petitot~\cite{Petitot_B}. Then this model was refined by G.~Citti and A.~Sarti~\cite{Citti_Sarti_B}. According to this model, the process of completion of occluded contours occurs by minimizing the excitation energy of neurons that perceive visual information in the areas of the observed scene that are occluded. This process can be interpreted as the action of the hypoelliptic diffusion operator studied in~\cite{Gauthier}. The resulting curves are sub-Riemannian length minimizers. Such curves are used for image inpainting~\cite{OptimalInpainting} and for explanation of some visual illusions~\cite{MashtakovDGA}.

In this article, we consider the problem of sub-Riemannian geodesics in the space $M$ of positions, orientations and curvatures. We formulate the problem as an optimal control problem in Section~\ref{sec:stat}. In Section~\ref{sec:control}, we prove complete controllability of the system and existence of optimal controls. Then, in Section~\ref{sec:PMP}, we apply a necessary optimality condition --- Pontryagin maximum principle (PMP), and examine the Hamiltonian system of PMP. In conclusion, we summarize the main results.

% Идея включения кривизны изначально была предложена Ж.~Петито~\cite{petitot}, где он предложил рассмотреть 4-мерное расширение группы $\SE(2)$

\section{Optimal control problem}\label{sec:stat}
We consider the following control system:
\begin{equation}\label{eq:contsyst}
%\begin{cases}
\left\{\begin{array}{l}
\dot{x} = u_1 \cos{\theta}, \\[0.005\linewidth]
\dot{y} = u_1 \sin{\theta}, \\[0.005\linewidth]
\dot{\theta} = u_1 k, \\[0.005\linewidth]
\dot{k} = u_2,
\end{array}\right.
%\end{cases}
\quad
\begin{array}{l}
(x,y,\theta, k) =q \in \SE(2) \times \R = M, \\[0.01\linewidth]
(u_1,u_2)  \in \R^2.
\end{array}
\end{equation}
For the trajectory corresponding to the control $(u_1(t), u_2(t))$, $t\in [0,T]$, $T>0$, we define a cost functional --- the sub-Riemannian length of this trajectory:
\begin{equation}\label{eq:int}
l = \int\limits_0^T \sqrt{u_1^2 (t)+u_2^2 (t) } \, dt.
\end{equation}
We study the problem of finding a Lipschitzian curve $\gamma:[0,T] \to M$ --- that is a trajectory of system~(\ref{eq:contsyst}), satisfying the given boundary conditions
\begin{equation}\label{eq:bounds}
\gamma(0)= q_0,\quad \gamma(T)=q_1, \qquad q_i \in M,
\end{equation}
and having the minimal sub-Riemannian length $l(\gamma)\to \min$.
\begin{remark}\label{rem:bounds}
It is easy to check that system~{\normalfont(\ref{eq:contsyst})} is invariant under parallel translations and rotations in the plane $(x,y)$. Because of this, without loss of generality, we can reduce the problem for an arbitrary $q_0=(x_0,y_0,\theta_0,k_0)$ to the case $q_0=(0,0,0,k_0).$
\end{remark}

\section{Complete controllability of the system}\label{sec:control}
The first question that arises when studying problem~(\ref{eq:contsyst})--(\ref{eq:bounds}) is existence of an admissible trajectory connecting boundary conditions~(\ref{eq:bounds}). If for any $q_0$, $q_1 \in M$ the answer is positive, then the control system is called completely controllable. Let us investigate controllability of system~(\ref{eq:contsyst}), using the technique of geometric control theory~\cite{notes}.

System~(\ref{eq:contsyst}) has the following form:
\begin{equation*}\label{eq:contsysvect}
\dot{\gamma} = u_1 X_1 +u_2 X_2,
\end{equation*}
where the vector fields near the controls are given by
\begin{equation*}\label{eq:vectfields}
%\begin{cases}
X_1 = 
\left(\begin{array}{c}
\cos{\theta} \\[0.005\linewidth]
\sin{\theta} \\[0.005\linewidth]
k \\[0.005\linewidth]
0\\
\end{array}\right),
\quad
%\end{cases}
%\begin{cases}
X_2 = 
\left(\begin{array}{c}
0 \\[0.005\linewidth]
0 \\[0.005\linewidth]
0 \\[0.005\linewidth]
1 \\
\end{array}\right).
%\end{cases}
\end{equation*}

We investigate controllability using Chow–Rashevskii theorem. In our case, it suffices to check that the rank condition is satisfied. To do this, we calculate the following Lie brackets of the fields $X_i$:
\begin{equation*}\label{eq:X3X4}
X_3=[X_1,X_2] = %0 - 
%\begin{pmatrix}
%0 & 0 & -\sin{\theta} & 0 \\
%0 & 0 & \cos{\theta} & 0\\
%0 & 0 & 0 & 1\\
%0 & 0 & 0 & 0\\
%\end{pmatrix}
%\left(\begin{array}{c}
%0 \\[0.005\linewidth]
%0 \\[0.005\linewidth]
%0 \\[0.005\linewidth]
%1 \\
%\end{array}\right) = 
\left(\begin{array}{c}
0 \\[0.005\linewidth]
0 \\[0.005\linewidth]
-1 \\[0.005\linewidth]
0 \\
\end{array}\right),\quad
%\end{equation*}
%\begin{equation*}\label{eq:X4}
X_4=[X_1,X_3]  = 
\left(\begin{array}{c}
-\sin{\theta} \\[0.005\linewidth]
\cos{\theta} \\[0.005\linewidth]
0 \\[0.005\linewidth]
0 \\
\end{array}\right).
\end{equation*}
 
 For the matrix composed of the vector fields $X_1, \ldots, X_4$ we have
\begin{equation}\label{eq:rank}
\det 
\begin{pmatrix}
\cos{\theta} & 0 & 0 & -\sin{\theta}\\
\sin{\theta} & 0 & 0 & \cos{\theta}\\
k & 0 & -1 & 0\\
0 & 1 & 0 & 0\\
\end{pmatrix}\equiv 1. 
\end{equation}
We conclude that the rank of the matrix is four, so the vector fields $X_i$ are linearly independent. So, they define a basis of the tangent space $ T_q M $ at every point $q$. Thus, wee see that all the conditions of the Chow–Rashevskii theorem are satisfied and we obtain the following result.
\begin{theorem}\label{thm:control}
Control system~{\normalfont(\ref{eq:contsyst})} is completely controllable.
\end{theorem}
\begin{remark}
Since $u_i$ are unbounded, condition~{\normalfont(\ref{eq:contsysvect})} is equivalent to
$\dot{\gamma} \in \Delta_\gamma = \spann(X_1(\gamma), X_2(\gamma)),$
where the family of planes $ \Delta $ is called the distribution. Due to condition~{\normalfont(\ref{eq:rank})}, the growth vector of the distribution~$\Delta$ equals $(2,3,4)$. Such systems are called structures of Engel type.
\end{remark}

Further, a question of existence of optimal trajectories arises: does there always exist an admissible trajectory satisfying conditions~(\ref{eq:bounds}), on which the functional~(\ref{eq:int}) reaches its minimum value? For our problem~(\ref{eq:contsyst})--(\ref{eq:bounds}) the answer is positive. Existence of optimal trajectories is guaranteed by Filippov theorem~\cite{notes,zelikin}.

\section{Pontryagin maximum principle}\label{sec:PMP}
Before proceeding to examination of extremal trajectories, let us reduce the problem under consideration to a simpler one. By virtue of Cauchy-Schwarz inequality, the original problem is equivalent to the problem of minimizing the action functional
\begin{equation}\label{eq:actionfunc}
J = \int\limits_0^T \frac{u_1^2 (t)+u_2^2 (t)}2 \, dt \to \min.
\end{equation}

Apply to problem~(\ref{eq:contsyst}),~(\ref{eq:bounds}),~(\ref{eq:actionfunc}) a necessary condition of optimality --- Pontryagin maximum principle (PMP). We introduce the Pontryagin function
\begin{equation*}\label{eq:pontHam}
h_{u}^\nu =\langle p, \sum\limits_{i=1}^2 u_i X_i\rangle + \nu/2 \sum\limits_{i=1}^2 u_i^2, \quad p \in T^*M, \,\, \nu \leq 0.
\end{equation*}
PMP states that if $(u(t), q(t))$, $t \in[0,T]$, is an optimal process, then the following conditions hold:
\begin{enumerate}
\item Hamiltonian system $\displaystyle \dot{p}=-\frac{\partial h_u^\nu}{\partial q},\,\,  \dot{q}=\frac{\partial h_u^\nu}{\partial p}$;
\item Maximum condition $h_{u(t)}^\nu (p(t),q(t)) ={\underset{u \in \R^2}{\max}}\, h_{u}^\nu (p(t),q(t)) $;
\item Nontriviality condition $(p(t),\nu)\neq(0,0)$ $\forall t \in [0,T]$.
\end{enumerate}
Denote $h_i=\langle p,X_i\rangle.$ The Pontryagin function takes the form
\begin{equation*}\label{eq:pontHam1}
h_{u}^\nu =u_1 h_1+u_2 h_2 + \nu/2 \left(u_1^2+u_2^2\right).
\end{equation*}

In formulation of PMP, without loss of generality, it suffices to consider two cases: $\nu=0$ --- abnormal case and $\nu=-1$ --- normal case. Next, we will consider both cases in detail.

\subsection{Abnormal case $\nu=0$.}\label{sec:abnormal}
The Pontryagin function is $h_{u}^0 =u_1 h_1+u_2 h_2.$ This is a linear function unbounded when $h_1^2+ h_2^2 \neq 0$. Thus, the maximum condition is satisfied if and only if $h_1=h_2\equiv0$. The maximized Hamiltonian in this case is $H = {\underset{u\in \R^2}{\max}} h_u^0 = 0.$

Let $a,b,c,d \in \R$ be the components of the covector $p$ in canonical Darboux coordinates. They are changing by the law $\dot{p_i}=-\frac{\partial h_{u}^0}{\partial q_i}$. 

The Hamiltonian system of PMP has the form
\begin{equation}\label{eq:vertpartabnormal}
\left\{\begin{array}{l}
\dot{x} = u_1 \cos\theta, \\[0.005\linewidth]
\dot{y} = u_1 \sin\theta, \\[0.005\linewidth]
\dot{\theta} = u_1 k, \\[0.005\linewidth]
\dot{k} = u_2,
\end{array}\right.
\quad
\left\{\begin{array}{l}
\dot{a} = 0, \\[0.005\linewidth]
\dot{b} = 0, \\[0.005\linewidth]
\dot{c} = -u_1 (-a \sin{\theta}+b \cos{\theta}), \\[0.005\linewidth]
\dot{d} = -u_1 \, c.
\end{array}\right.
\end{equation}

The condition $h_1=h_2\equiv0$ implies
\begin{equation}\label{eq:h1h2eq0}
\left\{\begin{array}{l}
a \cos{\theta} + b \sin{\theta} + c k \equiv0, \\[0.005\linewidth]
d \equiv 0.
\end{array}\right.
\end{equation}

%В силу Замечания~\ref{rem:bounds} первое тождество этой системы должно выполняться в частности при $\theta=0$, то есть $a+c\, k \equiv 0$. При этом третье уравнение системы~(\ref{eq:vertpartabnormal}) принимает вид $\dot{c}= -u_1 b$.

By virtue of system~(\ref{eq:vertpartabnormal}), the second identity implies $\dot{d} = -u_1 \, c \equiv 0$, while differentiation of the first identity implies $u_2 \, c \equiv 0$. Thus, we have 
\begin{equation}\label{eq:abnormalu}
(u_1^2+u_2^2) c^2 \equiv 0.
\end{equation} %\Leftrightarrow \forall t\in[0,T]: \left[\begin{array}{l} u_1(t)=u_2(t)=0,\\ c(t)=0.  \end{array} \right. $$ 

Since every admissible curve of positive length is a Lipschitz reparameterization of an arclength parameterized admissible one, see \cite[Lemma 3.16]{b:agrachevbook}, one can chose the natural parameterization $u_1^2+u_2^2=1$ for every trajectory of system~(\ref{eq:contsyst}) that is not a fixed point. Thus, identity~(\ref{eq:abnormalu}) is equivalent to $c\equiv0$ on the intervals of time where the trajectory is not a fixed point.

The identity $c \equiv 0$ implies $\dot{c}\equiv0$. From the third equation of system~(\ref{eq:vertpartabnormal}) it follows $u_1 (-a \sin{\theta}+b \cos{\theta}) \equiv 0$. On the other hand, the first identity of system~(\ref{eq:h1h2eq0}) takes the form $a \cos{\theta} + b \sin{\theta} \equiv0$. Thus, we have
$u_1^2(a^2+b^2)\equiv0$. Note, that if $a^2+b^2 = 0$ then  $a=b=c=d=0$, which contradicts the nontriviality condition in PMP. Thus, it remains to consider the case $u_1 \equiv 0$. 

%This identity holds when either $u_1 \equiv 0$, which leads to the case considered below, or $(-a \sin{\theta}+b \cos{\theta}) \equiv 0$.  We conclude that $a^2+b^2 \equiv 0$. Thus, we have $a=b=c=d\equiv0$, which contradicts the nontriviality condition in PMP. 

%Since $x$, $y$, $\theta$, $k$, $a$, $b$, $c$, $d$ are Lipschitzian functions they are differentiable for almost every $t \in [0,T]$.  

%Two cases are possible: $c \equiv 0$ or $u_1 \equiv 0$.

%Let's consider the case $c \equiv 0$. If a function is identically equal to zero, then its derivative is also equal to zero. Thus we have $\dot{c}\equiv0$. From the third equation of system~(\ref{eq:vertpartabnormal}) it follows $u_1 (-a \sin{\theta}+b \cos{\theta}) \equiv 0$. This identity holds when either $u_1 \equiv 0$, which leads to the case considered below, or $(-a \sin{\theta}+b \cos{\theta}) \equiv 0$. On the other hand, the first identity of system~(\ref{eq:h1h2eq0}) takes the form $a \cos{\theta} + b \sin{\theta} \equiv0$. We conclude that $a^2+b^2 \equiv 0$. Thus, we have $a=b=c=d\equiv0$, which contradicts the nontriviality condition in PMP. 

In this case, the covector $ p $ is constant and nonzero. We obtain that the abnormal extremals have the form $\gamma(t)=(0,0,0,k_0+U(t))$, where $U(t)= \int_0^t u_2(\tau) d\tau$, and $u_2(t)$ is any real-valued $L^\infty(0,T)$ function. It is easy to see that the trajectory is not optimal, if $u_2$ changes its sign. This holds since when $u_2$ changes its sign, the trajectory is followed in opposite directions. Choosing the natural parameterization $u_1^2+u_2^2=1$ on optimal trajectories, we obtain $u_2=\pm 1$, that results to the optimal trajectory $\gamma(t)=(0,0,0,k_0\pm t)$. Thus, we obtain the following result.

\begin{theorem}\label{thm:abnormal}
Abnormal extremal trajectories in problem~{\normalfont(\ref{eq:contsyst})--(\ref{eq:bounds})} have the form $\gamma(t) = (0,0,0,k_0+\int_0^t u_2(\tau) d \tau)$, where $u_2(t)$ is any real-valued $L^\infty(0,T)$ function. Naturally parameterized abnormal optimal trajectories have the form $\gamma(t) = (0,0,0,k_0\pm t)$. 
\end{theorem}

\subsection{Normal case $\nu=-1$.}\label{sec:normal}
The Pontryagin function $h_{u}^{-1}$ takes the form $h_{u}^{-1} =u_1 h_1+u_2 h_2 - \left( u_1^2+u_2^2\right)/2.$  The maximum condition gives the expression for the extremal controls:
$$\frac{\partial h_{u}^{-1}}{\partial u_i} = h_i-u_i=0 \, \Rightarrow \, u_i=h_i.$$ 
The maximized Hamiltonian $H = {\underset{u\in \R^2}{\max}}h_{u}^{-1}$ takes the form $$H= \frac{h_1^2+h_2^2}{2}.$$

By definition $h_i = \langle p, X_i \rangle$, $i=1,\ldots,4$ we have
\begin{equation*}\label{eq:hdef}
%\begin{cases}
\left\{\begin{array}{l}
h_1 = a \cos{\theta} + b \sin{\theta} + c k, \\[0.005\linewidth]
h_2 = d, \\[0.005\linewidth]
h_3 = -c, \\[0.005\linewidth]
h_4= -a \sin{\theta} + b \cos{\theta}.
%,\\
%\dot{a} = 0, \\[0.005\linewidth]
%\dot{b} = 0, \\[0.005\linewidth]
%\dot{c} = -((b^2-a^2)\cos{\theta}\sin{\theta}+ab(\cos^2{\theta-\sin^2{\theta}}) -ack\sin{\theta}+bck\cos{\theta}), \\[0.005\linewidth]
%\dot{d} = -(c^2k+ac\cos{\theta}+bc\sin{\theta}),\\
%\dot{x} = (a \cos{\theta} + b \sin{\theta} + ck)\cos{\theta}, \\[0.005\linewidth]
%\dot{y} = (a \cos{\theta} + b \sin{\theta} + ck)\sin{\theta}, \\[0.005\linewidth]
%\dot{\theta} = (a \cos{\theta} + b \sin{\theta} + ck)k, \\[0.005\linewidth]
%\dot{k} = d.
\end{array}\right.
\, \Leftrightarrow \,
\left\{\begin{array}{l}
a = (h_1 + k h_3) \cos\theta - h_4 \sin\theta,  \\[0.005\linewidth]
b = h_4 \cos\theta + (h_1 + k h_3) \sin\theta, \\[0.005\linewidth]
c = -h_3, \\[0.005\linewidth]
d =  h_2.
\end{array}\right.
%\end{cases}
\end{equation*}

The Hamiltonian system in canonical coordinates has the form
\begin{equation*}\label{eq:hamsysnormcanon}
\left\{\begin{array}{l}
\dot{x} = (a \cos{\theta} + b \sin{\theta} + c k)\cos{\theta}, \\[0.005\linewidth]
\dot{y} = (a \cos{\theta} + b \sin{\theta} + c k)\sin{\theta}, \\[0.005\linewidth]
\dot{\theta} = k (a \cos{\theta} + b \sin{\theta} + c k), \\[0.005\linewidth]
\dot{k} = d,
\end{array}\right.
\quad
\left\{\begin{array}{l}
\dot{a} = 0, \\[0.005\linewidth]
\dot{b} = 0, \\[0.005\linewidth]
\dot{c} =-(b \cos\theta - a \sin\theta)(k c + a \cos\theta+b \sin\theta), \\[0.005\linewidth]
\dot{d} = -c (c k+a \cos{\theta}+b \sin{\theta}).\\
\end{array}\right.
\end{equation*}

Rewrite this system in the coordinates $h_i$:
\begin{equation}\label{eq:hamsysnormh}
\left\{\begin{array}{l}
\dot{x} = h_1\cos{\theta}, \\[0.005\linewidth]
\dot{y} = h_1\sin{\theta}, \\[0.005\linewidth]
\dot{\theta} = h_1 k,\\[0.005\linewidth]
\dot{k} = h_2,
\end{array}\right.
\quad
\left\{\begin{array}{l}
\dot{h}_1 = -h_2 h_3, \\[0.005\linewidth]
\dot{h}_2 = h_1 h_3, \\[0.005\linewidth]
\dot{h}_3 = h_1 h_4, \\[0.005\linewidth]
\dot{h}_4= -k h_1 (k h_3+h_1).
\end{array}\right.
\end{equation}
By choosing the natural parameterization $u_1^2+u_2^2 =1$ on the extremal trajectories one fixes the level surface of the Hamiltonian $H = (h_1^2+h_2^2)/2 = 1/2$. Introduce in the plane $(h_1,h_2)$ a polar angle $\alpha \in S^1$:
$$h_1 = \cos\alpha, \quad h_2 = \sin\alpha.$$
By rewriting Hamiltonian system~(\ref{eq:hamsysnormh}) we get the following result.
\begin{theorem}\label{thm:normal}
Naturally parameterized normal extremal trajectories in problem~{\normalfont(\ref{eq:contsyst})--(\ref{eq:bounds})} are solutions to the system
\begin{equation}\label{eq:hamsysnorm}
\left\{\begin{array}{l}
\dot{x} = \cos\alpha\cos{\theta}, \\[0.005\linewidth]
\dot{y} = \cos\alpha\sin{\theta}, \\[0.005\linewidth]
\dot{\theta} = k \cos\alpha,\\[0.005\linewidth]
\dot{k} = \sin\alpha,
\end{array}\right.
\quad
\left\{\begin{array}{l}
\dot{\alpha} = h_3, \\[0.005\linewidth]
\dot{h}_3 = h_4 \cos\alpha, \\[0.005\linewidth]
\dot{h}_4= -k \cos\alpha(k h_3+\cos\alpha).
\end{array}\right.
\end{equation}
\end{theorem}

The question of existence of analytic expression for extremal trajectories arises: is the Hamiltonian system integrable? To prove Liouville integrability of system~(\ref{eq:hamsysnormh}) it suffices to find four functionally independent first integrals in involution. One such integral is the Hamiltonian $ H $, two more first integrals $ a $ and $ b $ follow directly from the representation of the Hamiltonian system in canonical coordinates. These three integrals are functionally independent and are in involution. The existence of the remainig first integral is questionable.

Numerical experiments have been carried out indicating the presence of the fourth integral. To this end, we consider the four-dimensional system, which is decoupled from the rest of the variables:
\begin{equation}\label{eq:foureq}
\left\{\begin{array}{l}
\dot{\alpha} = h_3, \\[0.005\linewidth]
\dot{h}_3 = h_4 \cos\alpha, \\[0.005\linewidth]
\dot{h}_4= -k \cos\alpha(k h_3+\cos\alpha), \\[0.005\linewidth]
\dot{k} = \sin\alpha.
\end{array}\right.
\end{equation}

Existence of the first integral of this system has been studied using Poincare map. For this, a periodic trajectory was found. Such a trajectory is given by the initial value $\alpha(0)=\frac{\pi}{2}$, $h_3(0)=1$, $h_4(0)=0$, $k(0)=0$. More precisely, a one-parameter family of periodic trajectories was found
\begin{equation*}\label{eq:loops}
\alpha(t)= \frac{\pi}{2}+t \, h_3(0),\,\, h_3(t) = h_3(0), \,\,  h_4(t) = 0, \, \,  k(t) = \frac{\sin(t \, h_3(0))}{h_3(0)},
\end{equation*}
 obtained by changing the initial value $h_3(0)$. In this case, the period is $T= |\frac{2 \pi}{h_3(0)}|$. With such initial values, the point  $(h_3,h_4) = (h_3(0),0)$ is fixed.

The hyperspace $k=0$, which is transversal to the periodic trajectory, was chosen. The Poincare map was constructed as follows. Initial values were chosen from a small neighborhood of initial values of the periodic trajectory. The trajectories departing from these initial values were computed. Points of intersections of every such trajectory with the transversal hyperspace were computed (for the first time, the second time, $\ldots$, the $N$-th time). 

 \begin{figure}[ht]
\centering
\includegraphics[width=0.49\linewidth]{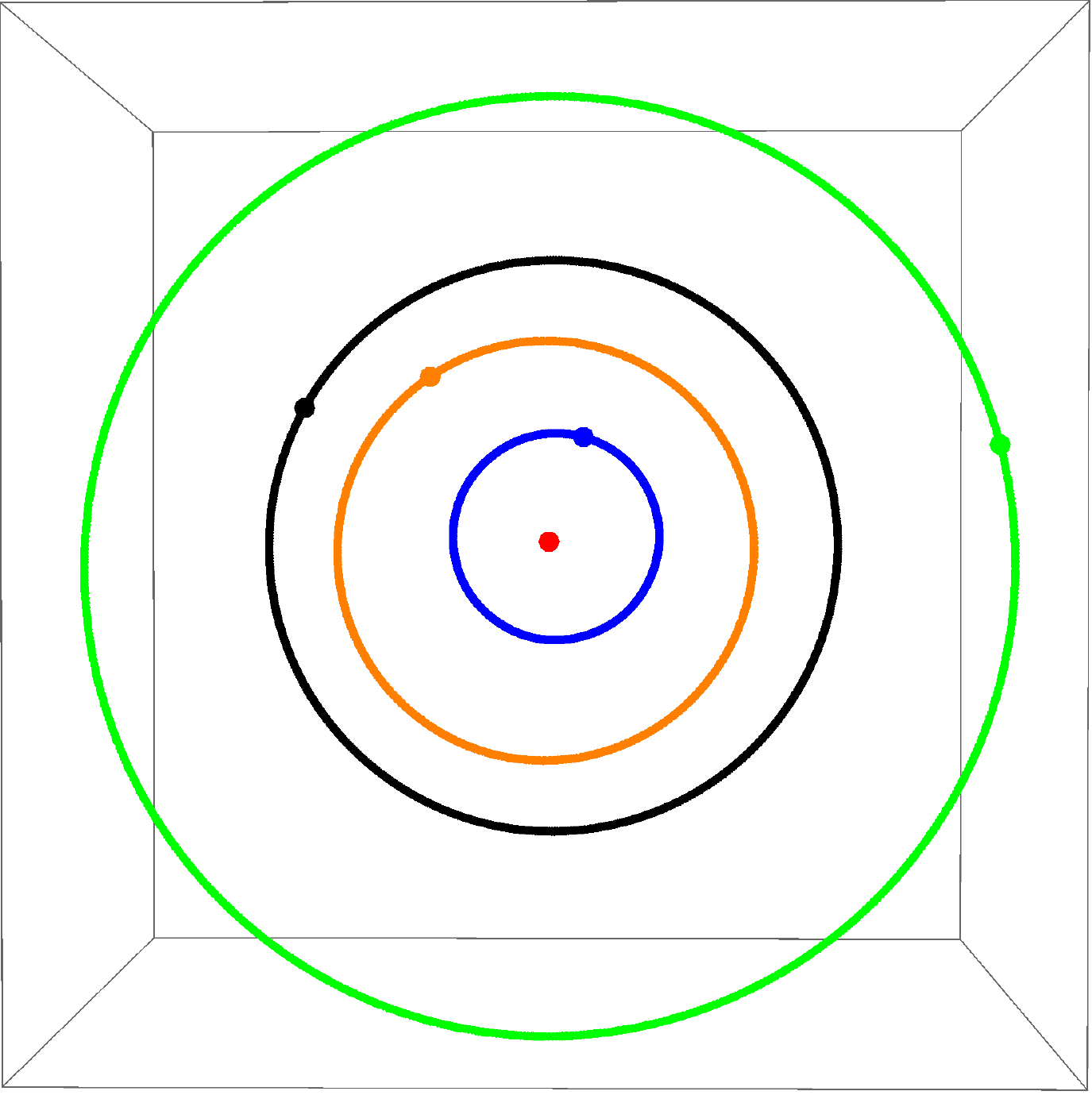}~
\includegraphics[width=0.49\linewidth]{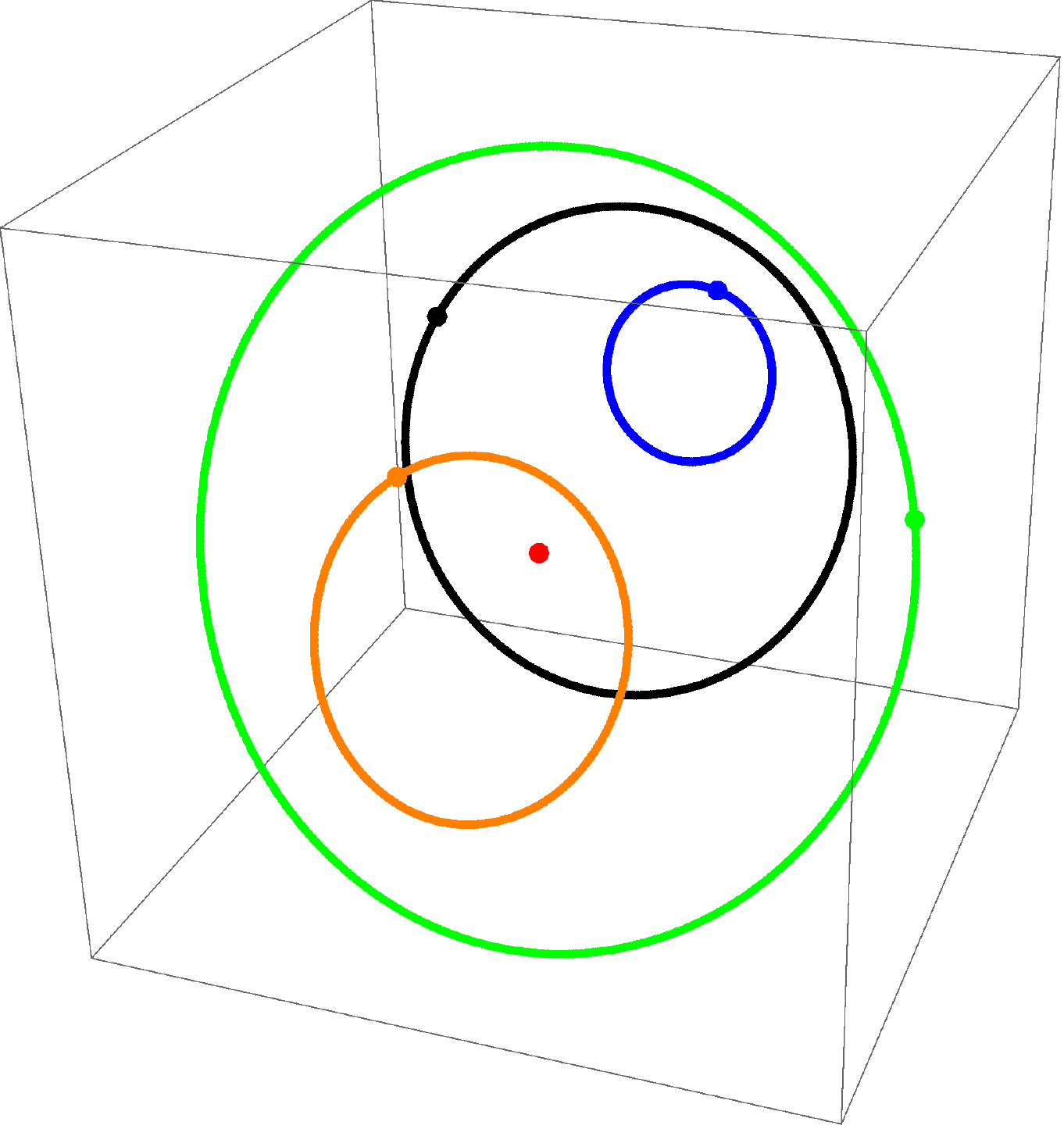}
\caption{Orbits of the Poincare map in the space $(\alpha,h_3,h_4)$ are formed by intersection points of the transversal hyperspace $k=0$ with trajectories close to the periodic one (red dot). Different orbits are depicted in different colors. Starting points are indicated for each trajectory.}
\label{fig:chaos}
\end{figure}

In Fig.~\ref{fig:chaos} the orbits of the Poincare map in the space $(\alpha,h_3,h_4)$ are shown. The red dot corresponds to the orbit of a periodic trajectory with initial value  $\alpha(0)=\frac{\pi}{2}$, $h_3(0)=1$, $h_4(0)=0$. For the rest of the trajectories, the following parameters were used. For the orange trajectory: $\alpha(0)=1.56$, $h_3(0) =0.94$, $h_4(0) =0.02$. For the green trajectory: $\alpha(0)=1.55$, $h_3(0) =1.06$, $h_4(0) =0.05$. For the black trajectory: $\alpha(0)=1.6$, $h_3(0) =1.14$, $h_4(0) =0.02$. For the blue trajectory: $\alpha(0)=1.58$, $h_3(0) =1.24$, $h_4(0) =0.01$. The number of iterations of the Poincare map for all trajectories was chosen $N=1000$.

The trajectories were computed by numerical integration (NDSolve in Wolfram Mathematica). The instances when the trajectory intersects the transversal hyperspace were determined using the numerical solution (FindRoot) of the equation $k(t)=0$ along the given trajectory with initial approximation $t= |\frac{2 \pi}{h_3(0)}|$ at each iteration $i=1,\ldots,N$. 

It is remarkable that the points of the Poincare map fill continuous closed curves. This indicates that for a small perturbation of the initial conditions, the trajectories remain close over a long time interval. This situation arises when the system is integrable (at least in some domain). %Since system~(\ref{eq:foureq}) is specified by analytic functions, then integrability in some domain implies integrability globally in the whole space. 

Thus, the presented numerical experiments indicate the presence of a fourth independent first integral, which results to Liouville integrability of the Hamiltonian system~(\ref{eq:hamsysnorm}).    

\section{Conclusion}\label{sec:concl}

In this paper, we consider the problem of sub-Riemannian geodesics in the four-dimensional manifold $M=\R^2\times\SO(2)\times\R$. This problem arises when modelling the mechanism of occluded contours completion in the extended model of the visual cortex by J.~Petitot, G.~Citti, and A.~Sarti. The extension of the classical 3D model is performed by taking the curvature of the contours into account.

The following main results are obtained in the article. Complete controllability of the system and existence of optimal controls are proved. The Hamiltonian system of PMP is derived. The explicit parametrization of abnormal trajectories is found. In the normal case, three functionally independent first integrals are found. A hypothesis is formulated, confirmed by numerical experiments, that the normal Hamiltonian system of PMP is Liouville integrable.

The authors are grateful to Prof. Yu.L. Sachkov for many fruitful suggestions and valuable comments on the work.

\end{document}